\documentclass[12pt]{article}
\usepackage[letterpaper,twoside,hcentering]{geometry}
\pdfoutput=1

\def\titlelong{Estimating a bivariate linear relationship}
\def\titleshort{Bivariate linear relationship}
\def\authorlong{David Leonard}
\def\authorshort{D.~Leonard}

\usepackage{fancyhdr}
\usepackage{mathptmx} 
\usepackage[scaled=0.90]{helvet} 
\usepackage{courier} 
\normalfont
\usepackage[T1]{fontenc}

\usepackage{amssymb}
\usepackage{amsmath}
\usepackage{graphicx}
\usepackage{url}
\usepackage{setspace}
\usepackage[authoryear,round]{natbib}

\DeclareMathOperator{\tr}{tr}

\DeclareMathOperator{\sign}{sign}
\DeclareMathOperator{\Nrml}{N}
\DeclareMathOperator{\W}{W}
\newcommand{\dx}[1]{\,\mathrm{d}{#1}}
\newcommand{\transpose}[1]{{#1}^\mathsf{T}}
\newcommand{\abs}[1]{\lvert{#1}\rvert}

\begin{document}

\title{
\vspace{-0.75in}\titlelong\footnote{
Originally published in {\em Bayesian Analysis} (2011) 6: 727--754, DOI:10.1214/11-BA627.
}}
\author{
\authorlong\footnote{
Department of Clinical Sciences, University of Texas Southwestern Medical Center, Dallas, TX, david.leonard@utsouthwestern.edu
}}
\date{}
\maketitle

\thispagestyle{empty}

\begin{abstract}
Solutions of the bivariate, linear errors-in-variables estimation problem with unspecified errors are expected to be invariant under interchange and scaling of the coordinates.  The appealing model of normally distributed true values and errors is unidentified without additional information.  I propose a prior density that incorporates the fact that the slope and variance parameters together determine the covariance matrix of the unobserved true values but is otherwise diffuse.  The marginal posterior density of the slope is invariant to interchange and scaling of the coordinates and depends on the data only through the sample correlation coefficient and ratio of standard deviations.  It covers the interval between the two ordinary least squares estimates but diminishes rapidly outside of it.  I introduce the R package \textsf{leiv} for computing the posterior density, and I apply it to examples in astronomy and method comparison.

\end{abstract}

\section{Introduction}\label{introduction}

Simple linear relationships inspire much empirical research, yet how to estimate their parameters is a topic of continuing debate.  Longstanding examples include the permanent income model in economics \citep{zellnerbook}, cosmic distance scale applications in astronomy \citep{isobe1990}, and allometric studies in biology \citep{warton2006}.  From a statistical perspective, the common goal of these investigations is to estimate the slope relating two variables that are observed with error.  The controversy stems from the absence of an estimate that is invariant to interchange and scaling of the coordinates and depends reasonably on their joint distribution.

In one of the earliest comprehensive reviews, \citet{madansky1959} fixes ideas with the familiar problem of estimating the density of a solid by fitting a line to measurements of the mass and volume of a number of specimens.  In this problem, the density estimate should not depend on which axes the variables are plotted.  It should also not depend on the units of measurement.  That is, the same inference should be made by applying a scale conversion to the data before fitting or to the density estimate afterward.  The observations may be affected by measurement errors as well as errors intrinsic to the specimens, such as contamination by unknown impurities.

In many applications, the linear relationship appears on the log-log scale.  In such instances, units of measurement do not affect the slope of the fitted line, so it may seem that scale invariance is unnecessary.  \citet{warton2006} point out, however, that many multiplicative relationships involve arbitrary powers of the variables that translate to scale changes upon log transformation.  They offer that in an allometric analysis of certain saplings, for example, analyzing the relationship between height and basal diameter or basal area should lead to the same scientific conclusions.  Similar considerations apply in the analysis of the Faber-Jackson relation, Section \ref{faberjackson}.

The ordinary least-squares (OLS) estimate is scale invariant but not invariant to interchange of the coordinates.  The orthogonal regression estimate, proposed by \citet{adcock1877,adcock1878} and \citet{pearson1901}, is invariant to interchange of the coordinates, but it was famously criticized by \citet{wald1940} for its lack of scale invariance.  The economist \citet{samuelson1942} proposed the additional property of dependence only on the sample correlation coefficient and ratio of standard deviations.  He showed that the only point estimate of the slope exhibiting these particular invariance and dependence properties is the geometric mean of the two OLS estimates, an estimate he credited to \citet{frisch1934}.  This estimate, which is equal to the ratio of standard deviations of the measurements, depends on their joint distribution only for its sign, however.

Dependence on the correlation coefficient and ratio of standard deviations is especially appealing in the model of normally distributed true values contaminated by normally distributed errors.  \citet{reiersol1950} demonstrated that this model is unidentified; the sampling density identifies not a point but a continuum of estimates.  In some situations, such as in pure measurement error problems, supplementing the data with replicate measurements may solve the identification problem.  In others, the additional information must come from outside the sample.  A prior density would provide a natural way to incorporate it.

This article will show how to assign a prior density jointly to the slope and variance parameters that leads to a marginal posterior density of the slope that is invariant under interchange and scaling of the coordinates and has sufficient statistics in the sample correlation coefficient and ratio of standard deviations.  Passage to an appropriate noninformative limit is possible at the very end of the calculation.  In contrast, previous Bayesian solutions have relied on independent, informative prior densities for the variance parameters \citep[see, for example,][]{zellnerbook,polasek1993}.

In one of the earliest reported Bayesian analyses, \citet{lindley1968} predicted some general properties of the marginal posterior density without fully specifying the prior density.  These properties, notably including failure to concentrate around a single value in the limit of infinitely large samples, are indeed exhibited by the fully specified solution that follows.

This introduction has mentioned only a few relevant developments in the long history of this problem.  \citet{madansky1959}, \citet{anderson1984}, \citet{sprent1990} and \citet{stefanski2000} provide comprehensive reviews.  The book by \cite{fuller1987} has become a standard reference.  \citet{isobe1990} describe applications in astronomy, biology, chemistry, geology and physics.

\section{Formulation of the problem}\label{problem}

I adopt the notation of Zellner's \citeyearpar[][Chapter 5]{zellnerbook} comprehensive presentation.  The data are $n$ pairs $\{ y_{1i},y_{2i} \}$, viewed as independent, noisy observations of their unobserved true values $\{ \xi_{1i},\xi_{2i} \}$
\begin{equation}\label{observations}
\begin{split}
y_{1i}&=\xi_{1i}+u_{1i},\\
y_{2i}&=\xi_{2i}+u_{2i},
\end{split}
\end{equation}
$i=1,\dotsc,n$.  The elements of each pair of true values are linearly related,
\begin{equation}\label{line}
\xi_{2i}=\alpha+\beta \xi_{1i},
\end{equation}
and I seek an estimate of the slope $\beta$ and possibly the intercept $\alpha$.

I consider the model in which the true values $\{ \xi_{1i} \}$ and $\{ \xi_{2i} \}$ are samples from a bivariate normal distribution, degenerate to the regression line, and the total errors $\{ u_{1i} \}$ and $\{ u_{2i} \}$ are independently and normally distributed with mean zero and respective variances $\sigma_1^2$ and $\sigma_2^2$.  Despite outward appearances, assuming that either $\{ \xi_{1i} \}$ or $\{ \xi_{2i} \}$ are samples from an improper constant density imposes severe restrictions on the resulting solutions.  \citet{zellnerbook} showed that the model in which $\{ \xi_{1i} \}$ are samples from an improper constant density produces the same estimates of $\beta$ and $\alpha$ as OLS regression of $y_2$ on $y_1$.  He explained that the infinite variance of the distribution of $\{ \xi_{1i} \}$ makes the variance of the distribution of errors $\{ u_{1i} \}$ negligible in comparison.  Likewise, the model in which $\{ \xi_{2i} \}$ are samples from an improper constant density produces the same estimates as OLS regression of $y_1$ on $y_2$.  Indeed, any model that assumes an improper constant density along a line in the plane of true values presupposes a direction of ignorable error.  The normal model may be appropriate in case such specific information is not available.

Denoting the distribution of $\{ \xi_{1i} \}$ by $\Nrml(\mu_1,\tau^2)$, the sampling distribution of the observations $\{ y_{1i},y_{2i} \}$ is bivariate normal with mean
\begin{equation}\label{mean}
\mu=(\mu_1,\alpha+\beta\mu_1)
\end{equation}
and covariance
\begin{equation}\label{covtotal}
\Sigma = 
\begin{pmatrix}
\tau^2+\sigma_1^2& \beta\tau^2\\
\beta\tau^2& \beta^2\tau^2+\sigma_2^2
\end{pmatrix}.
\end{equation}

\citet{reiersol1950} demonstrated the consequence of the fact that the sampling distribution has six unknown parameters, but the sample mean and covariance matrix provide only five sufficient statistics.  Put simply, $\beta$ is not identifiable in the normal model without additional information.  OLS regression overcomes the difficulty by assuming one of $\sigma_1$ or $\sigma_2$ is zero, while orthogonal regression assumes the ratio $\sigma_2/\sigma_1=1$ or, more generally, a known constant.  These assumptions are more than adequate to estimate $\beta$; reducing the number of unknown parameters by one provides estimates of the remaining five.  The focus of the present effort is to find out how much the form of $\Sigma$ can tell us about $\beta$ alone.

\section{The posterior probability density}\label{posterior}

Ultimately, I will estimate $\beta$ from the posterior density
\begin{equation}\label{bayes}
p(\beta\mid y)=\frac{p(y\mid\beta)p(\beta)}{\int_{-\infty}^\infty p(y\mid\beta)p(\beta)\dx{\beta}},
\end{equation}
where $y$ is the $n \times 2$ observation matrix $(y_1,y_2)$, $p(\beta)$ is a prior density defined later in this section, and
\begin{equation}\label{sampling}
\begin{split}
p(y\mid\beta)&=\idotsint p(y\mid\mu_1,\alpha,\beta,\tau^2,\sigma_1^2,\sigma_2^2)\\
&\times p(\mu_1,\alpha,\tau^2,\sigma_1^2,\sigma_2^2\mid\beta) \dx{\mu_1}\dx{\alpha}\dx{\tau^2}\dx{\sigma_1^2}\dx{\sigma_2^2}
\end{split}
\end{equation}
is the reduced sampling density.

As discussed in Section~\ref{problem}, the full sampling density in the integrand of \eqref{sampling} is bivariate normal
\begin{equation}\label{normal}
\begin{split}
p(y&\mid\mu_1,\alpha,\beta,\tau^2,\sigma_1^2,\sigma_2^2)\\
&=(2\pi)^{-n/2}\abs{\Sigma}^{-n/2}\exp{\{-\tfrac{1}{2}\tr[(n\transpose{(\bar{y}-\mu)}(\bar{y}-\mu)+\nu S)\Sigma^{-1}]\}},
\end{split}
\end{equation}
where $\bar{y}=(\bar{y}_1,\bar{y}_2)$ is the vector of sample means, $\nu=n-1$, and $S$ is the sample covariance matrix with divisor $n-1$.

I factor the conditional prior density of the location and variance parameters in \eqref{sampling} as
\begin{equation}\label{prior}
p(\mu_1,\alpha,\tau^2,\sigma_1^2,\sigma_2^2\mid\beta) = p(\mu_1,\alpha\mid\tau^2,\sigma_1^2,\sigma_2^2,\beta) p(\tau^2,\sigma_1^2,\sigma_2^2\mid\beta)
\end{equation}
and take the conditional prior density of the location parameters $\mu_1$ and $\alpha$ to be a constant.

Previous Bayesian analyses have gone forward under the assumption that various functions of the variance parameters, for example, the ratio $\sigma_2/\sigma_1$, are approximately known (see, for example, \citealp[][Section~5.4]{zellnerbook}; \citealp{polasek1993}).  In the absence of such knowledge, the fact that $\Sigma$ links the variance parameters to the slope should not be ignored.  In particular, from \eqref{covtotal}, nonnegativity of the variance parameters
\begin{equation}\label{nonnegativity}
\begin{split}
\tau^2 &= \frac{1}{\beta}\Sigma_{12} \ge 0,\\
\sigma_1^2 &= \Sigma_{11}-\frac{1}{\beta}\Sigma_{12} \ge 0,\\
\sigma_2^2 &= \Sigma_{22}-\beta\Sigma_{12} \ge 0,
\end{split}
\end{equation}
modifies the domain of $\Sigma$ given $\beta$, and this information can be incorporated by assigning a conditional prior density to $\Sigma$ given $\beta$ after changing variables $p(\tau^2,\sigma_1^2,\sigma_2^2\mid\beta)=\abs{\beta}p(\Sigma\mid\beta)$.

Assigning an inverted Wishart density to $p(\Sigma\mid\beta)$ acknowledges that $\Sigma$ generates the Wishart distributed sample covariance $S$ \citep[][Chapter~7]{andersonbook}.  However, because the domain of $\Sigma$ depends on $\beta$, the inverted Wishart form
\begin{equation}\label{iwishart}
p(\Sigma\mid\beta)=K(\beta,\nu_0,\Psi_0(\beta))^{-1}\abs{\Sigma}^{-(\nu_0+3)/2}\exp{[-\tfrac{1}{2}\tr(\Psi_0(\beta)\Sigma^{-1})]}
\end{equation}
features a normalization factor that is a function of $\beta$.  It also introduces a degrees of freedom parameter $\nu_0$, a correlation parameter $\rho_0$, and a scale parameter $\kappa_0$ with units of $y_1$ through the precision matrix 
\begin{equation}\label{precision}
\Psi_0(\beta)=\nu_0 \kappa_0^2
\begin{pmatrix}
1& \rho_0\beta\\
\rho_0\beta& \beta^2
\end{pmatrix}.
\end{equation}
Fortunately, it will be possible to take the limit $\nu_0 \rightarrow 0$ at the very end of the calculation, removing any information these parameters carry, for any $-1 < \rho_0 < 1$ and $\kappa_0 > 0$.  The normalization factor $K(\beta,\nu_0,\Psi_0(\beta))$ is the crux of the method; it is calculated in Appendix~1.

After using \eqref{normal}--\eqref{iwishart} to carry out the integrations in \eqref{sampling}, the posterior density \eqref{bayes} is
\begin{equation}\label{bayes1}
p(\beta\mid y)=\frac{p(\beta)K(\beta,\nu_0+\nu,\Psi_0(\beta)+\nu S)/K(\beta,\nu_0,\Psi_0(\beta))}{\int_{-\infty}^\infty p(\beta)K(\beta,\nu_0+\nu,\Psi_0(\beta)+\nu S)/K(\beta,\nu_0,\Psi_0(\beta))\dx{\beta}}.
\end{equation}
Appendix~2 shows that in the limit $\nu_0 \rightarrow 0$,
\begin{equation}\label{posteriorlimit}
p(\beta\mid y)=\frac{p(\beta)J(\beta,\nu,r,l)}{\int_{-\infty}^\infty p(\beta)J(\beta,\nu,r,l)\dx{\beta}},
\end{equation}
where the sample correlation coefficient $r=S_{12}/(S_{11} S_{22})^{1/2}$ and the ratio of standard deviations $l=(S_{22}/S_{11})^{1/2}$ are sufficient statistics, and
\begin{equation}\label{Jbody}
J(\beta,\nu,r,l)=I(\abs{\beta}/l,\nu,r\sign(\beta))+I(l/\abs{\beta},\nu,r\sign(\beta)).
\end{equation}
In \eqref{Jbody},
\begin{equation}\label{Ibody}
I({\tilde\beta},\nu,r)=\int_{t_-(\nu,r)}^{t_+({\tilde\beta},\nu,r)} p_t(t;\nu)P_F(F(t,{\tilde\beta},\nu,r);\nu+1,\nu-1)\dx{t},
\end{equation}
where $p_t(t;\nu)$ is the Student $t$ probability density function, $P_F(F;\nu_1,\nu_2)$ is the $F$ cumulative distribution function,
\begin{align}
t_-(\nu,r)&=-\sqrt{\nu}\,r/\sqrt{1-r^2}, \label{tnbody}\\
t_+({\tilde\beta},\nu,r)&=\sqrt{\nu}\,({\tilde\beta}-r)/\sqrt{1-r^2}, \label{tpbody}
\end{align}
and
\begin{equation}\label{Fbody}
F(t,{\tilde\beta},\nu,r)=\frac{\nu-1}{\nu+1}\,\frac{\nu+t^2}{[t_+({\tilde\beta},\nu,r)-t_-(\nu,r)]^2-[ t-t_-(\nu,r)]^2}.
\end{equation}

Importantly, \eqref{Jbody} shows that $J(\beta,\nu,r,l)=J(\beta/l,\nu,r,1)$, so the entire role of $\beta$ is mediated by the scale invariant parameter ${\tilde\beta}=\beta/l$.  In terms of ${\tilde\beta}$ the posterior density \eqref{posteriorlimit} is
\begin{equation}\label{posteriorbeta}
p(\beta\mid y)=p({\tilde\beta}\mid y)/l,
\end{equation}
where
\begin{equation}\label{posteriorbetatilde}
p({\tilde\beta}\mid y)=\frac{p({\tilde\beta})J({\tilde\beta},\nu,r,1)}{\int_{-\infty}^\infty p({\tilde\beta})J({\tilde\beta},\nu,r,1)\dx{\tilde\beta}}
\end{equation}
depends on the data only through the sample correlation coefficient.

It remains to assign the prior density $p({\tilde\beta})$.  The pure number ${\tilde\beta}$ ensures scale invariance.  At a minimum, the prior specification should be invariant under interchange of the coordinates.  The sampling density \eqref{normal}, however, is invariant under continuous rotations of the coordinate plane, and for now I assume that the prior information available on ${\tilde\beta}$ is indifferent to such rotations as well.  I will return to this point briefly in Section~\ref{discussion}.   A rotationally invariant prior density will necessarily be invariant under interchange of the coordinates, by rotation through the angle $\pi/2$.  

Under rotation of the coordinates through an angle $\varphi$, a rotationally invariant prior density must satisfy the functional equation $p({\tilde\beta})=p({\tilde\beta}^{\prime})\abs{\mathrm{d}{\tilde\beta}^{\prime}/\mathrm{d}{\tilde\beta}}$, where ${\tilde\beta}^{\prime}=({\tilde\beta}\cos\varphi-\sin\varphi)/(\cos\varphi+{\tilde\beta}\sin\varphi)$.  Conveniently, there is only one solution, the Cauchy density
\begin{equation}\label{cauchy}
p({\tilde\beta})=\frac{1}{\pi}\,\frac{1}{1+{\tilde\beta}^2},
\end{equation}
equivalent to a uniform density on the angle $\theta=\arctan{\tilde\beta}$.  Appendix~3 shows that the resulting posterior density \eqref{posteriorbeta}--\eqref{cauchy} has precisely the form required for a density that depends on the data only through the sample correlation coefficient and ratio of standard deviations to be invariant under interchange and scaling of the coordinates.

The function $J(\beta,\nu,r,l)$ in \eqref{Jbody} is a sum of two integrals, one over the sampling density of the estimate $r l$ of $\beta$ in the OLS regression of $y_2$ on $y_1$ with $\nu$ degrees of freedom, and the other over the sampling density of the estimate $r/l$ of $1/\beta$ in the OLS regression of $y_1$ on $y_2$ with $\nu$ degrees of freedom.  These integrals are well-defined for $n>2$.  As $n$ becomes large, the Student $t$ density in the integrand of $I(\beta/l,\nu,r)$ becomes more sharply peaked around $t=0$, while the $F$ cumulative distribution function becomes more like a unit step function at $F(t,\beta/l,\nu,r)=1$.  Consequently, this integral contributes little to $J(\beta,\nu,r,l)$ unless the point $t=0$ is in $\{ t: F(t,\beta/l,\nu,r) \geq 1 \}$.  From definitions \eqref{tnbody}--\eqref{Fbody}, this condition is met whenever $\abs{\beta} \leq l$.  By the same reasoning, the integral $I(l/\beta,\nu,r)$ contributes little to $J(\beta,\nu,r,l)$ unless $\abs{\beta} \geq l$.  In other words, for $\abs{\beta}<l$, the posterior density is based largely on the sampling density of the estimate $r l$ of $\beta$ in the OLS regression of $y_2$ on $y_1$, whereas for $\abs{\beta}>l$, it is based largely on the sampling density of the estimate $r/l$ of $1/\beta$ in the OLS regression of $y_1$ on $y_2$.

In special cases, the integrals in $J(\beta,\nu,r,l)$ can be evaluated analytically.  For instance, Appendix~4 provides closed form expressions for the posterior density \eqref{posteriorbeta} for sample sizes of $n=4$ and $n=6$.  For $n=4$, the posterior density of $\theta=\arctan{\tilde\beta}$ is proportional to $\abs{\sin{2\theta}}/(1-r\sin{2\theta})$, $-\pi/2 \le \theta \le \pi/2$.  For $n=6$, the posterior density of $\theta$ is proportional to $\abs{\sin{2\theta}}(2-r\sin{2\theta})/(1-r\sin{2\theta})^2$.  These densities have relative maxima at $\theta=\pm \pi/4$ and absolute maximum at $\theta=\sign(r)\,\pi/4$.

\begin{figure}[!ht]
\begin{center}
\includegraphics[width=0.85\textwidth]{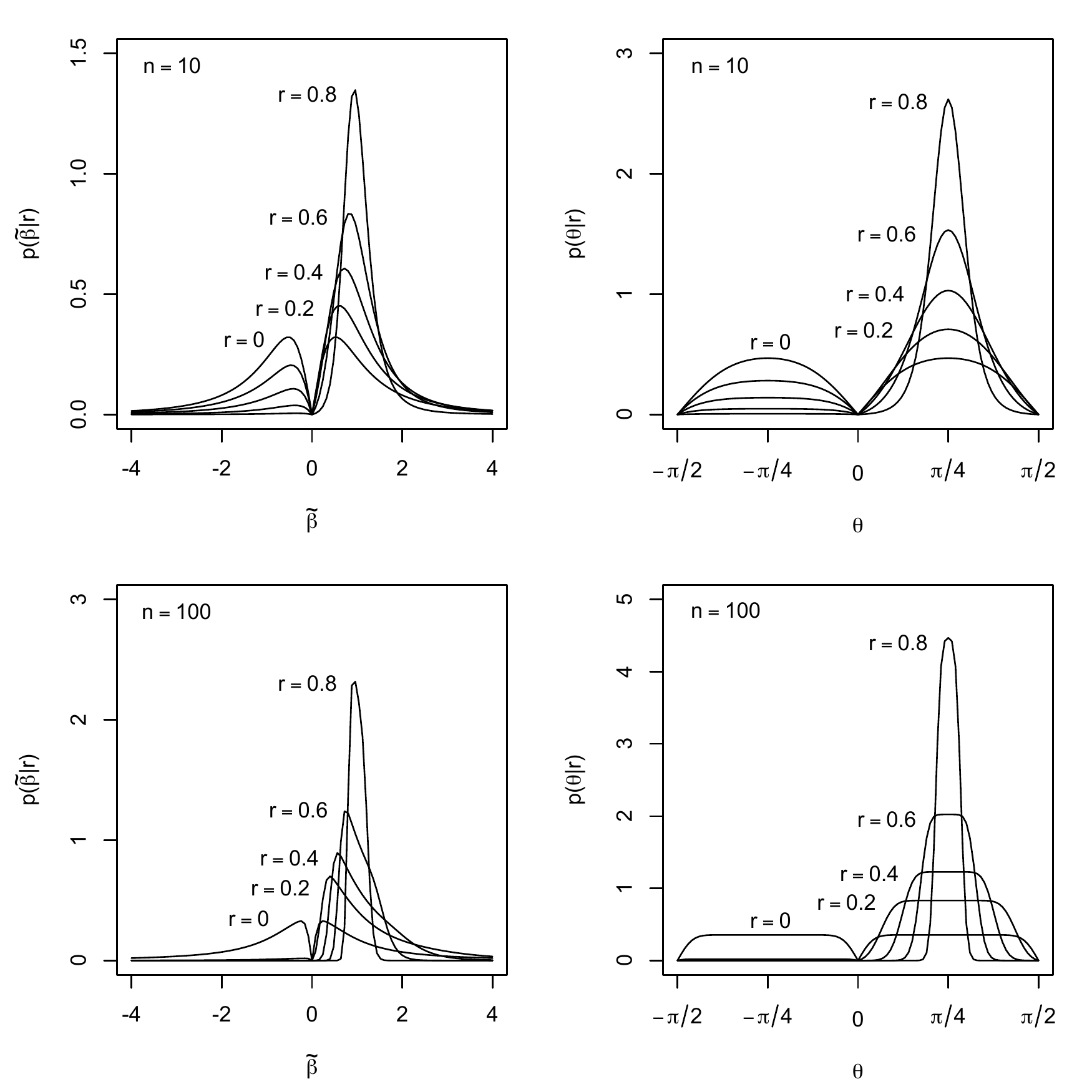}
\caption{Posterior densities of the scale invariant slope ${\tilde\beta}$ (left) computed from \eqref{posteriorbetatilde}--\eqref{cauchy} and of the corresponding angle $\theta=\arctan{\tilde\beta}$ (right) for sample sizes of $n=10$ (upper) and $n=100$ (lower) and a series of sample correlation coefficients.}\label{densityplots}
\end{center}
\end{figure}

More generally, the posterior density of the scale invariant slope \eqref{posteriorbetatilde} is illustrated in Figure~\ref{densityplots} for sample sizes of $n=10$ and $n=100$.  Notable features include the symmetry about ${\tilde\beta}=0$ for $r=0$ and the concentration about ${\tilde\beta}=\pm 1$ as $r \rightarrow \pm 1$.  For $\abs{r}<1$, however, the width does not shrink to zero as $n \rightarrow \infty$.  Figure~\ref{densityplots} also shows the posterior density of the corresponding angle $\theta=\arctan{\tilde\beta}$, for which the prior density \eqref{cauchy} is uniform.  The R package \textsf{leiv} \citep{R,leiv} computes the posterior density \eqref{posteriorbeta}--\eqref{cauchy} and is freely available from the Comprehensive R Archive Network (CRAN).

\section{Examples and simulations}\label{examples}

\subsection{Zellner's artificial data}\label{zellner}

\begin{figure}[!ht]
\begin{center}
\includegraphics[width=0.85\textwidth]{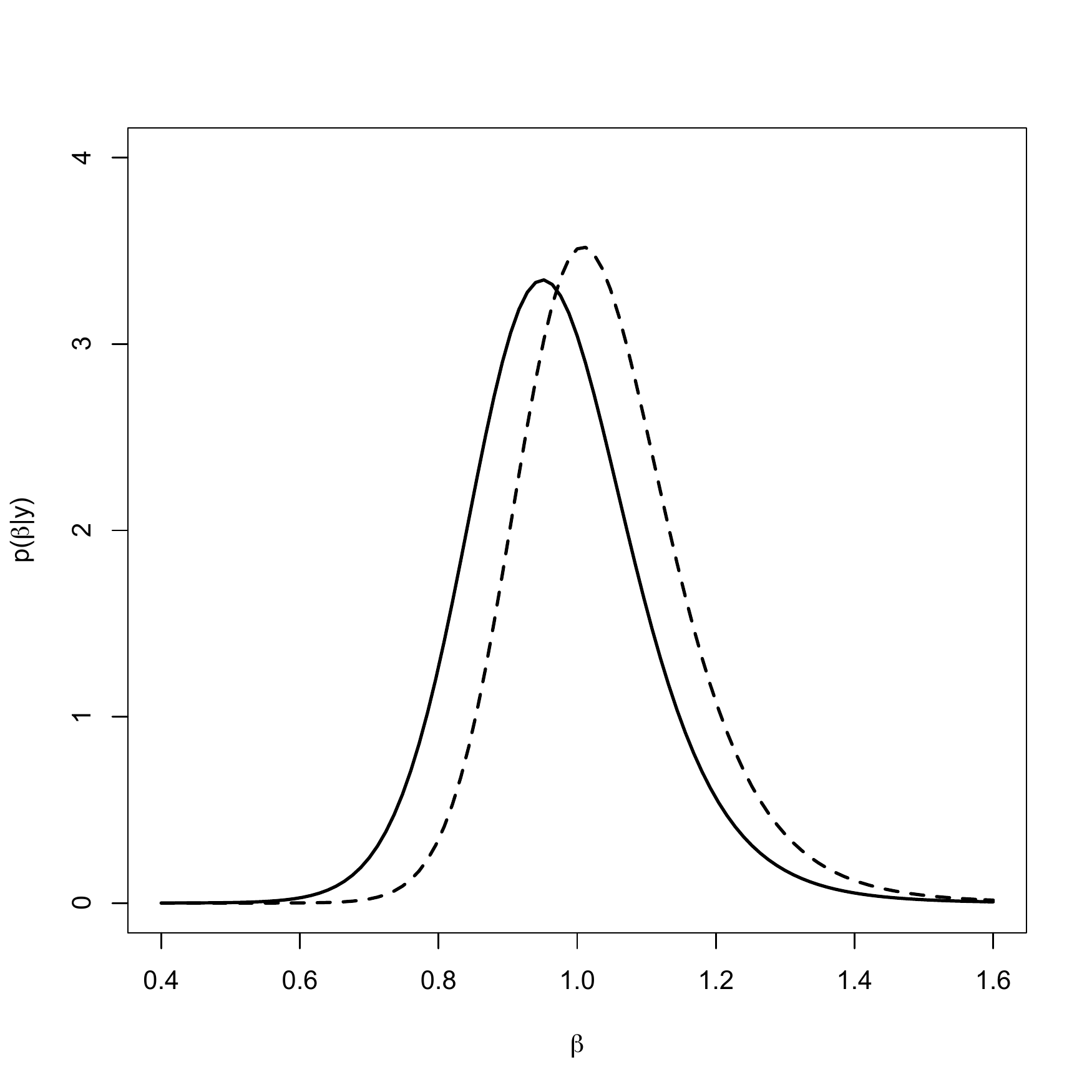}
\caption{Posterior density of the slope $\beta$ (solid) calculated from \eqref{posteriorbeta}--\eqref{cauchy} and (dashed) calculated by \citet[Figure~5.4]{zellnerbook} using an informative prior density.}\label{pZellner}
\end{center}
\end{figure}

The present example, using the artificial data of \citet[Table~5.1]{zellnerbook}, compares the posterior density \eqref{posteriorbeta} to Zellner's informed solution.  The data are $n=20$ pairs $\{ y_{1i},y_{2i} \}$, generated from the model \eqref{observations} and \eqref{line}, with slope $\beta=1$, intercept $\alpha=2$, error variances $\sigma_1^2=4$ and $\sigma_2^2=1$, true means $\mu_1=5$ and $\mu_2=\alpha+\beta\mu_1=7$, and true variance $\tau^2$=16.  These data meet all the assumptions of Section~\ref{problem}.  The sufficient statistics are $r=\text{0$\cdot$909}$ and $l=\text{0$\cdot$963}$.  The posterior density \eqref{posteriorbeta} is plotted in Figure~\ref{pZellner}.  The posterior median is 0$\cdot$963; the shortest 95\% probability interval is $(\text{0$\cdot$722},\text{1$\cdot$237})$.  For comparison, the 95\% confidence intervals are $(\text{0$\cdot$676},\text{1$\cdot$075})$ from the OLS regression of $y_2$ on $y_1$ and $(\text{0$\cdot$864},\text{1$\cdot$372})$ from the OLS regression of $y_1$ on $y_2$.

Figure~\ref{pZellner} also illustrates the posterior density that \citet[Figure~5.4]{zellnerbook} calculated from the same data, assuming a uniform prior density for $\beta$ and an independent, inverted gamma prior density for the true error variance ratio with mean 0$\cdot$246 and standard deviation 0$\cdot$152.  Zellner's posterior density is slightly narrower, due primarily to the informative prior density on the variance ratio.  It is shifted somewhat to the right, due in part to the uniform prior density for $\beta$, which is not rotationally invariant and favors angles approaching $\pm \pi/2$.

\subsection{Faber-Jackson relation}\label{faberjackson}

\begin{figure}[!ht]
\begin{center}
\includegraphics[width=0.85\textwidth]{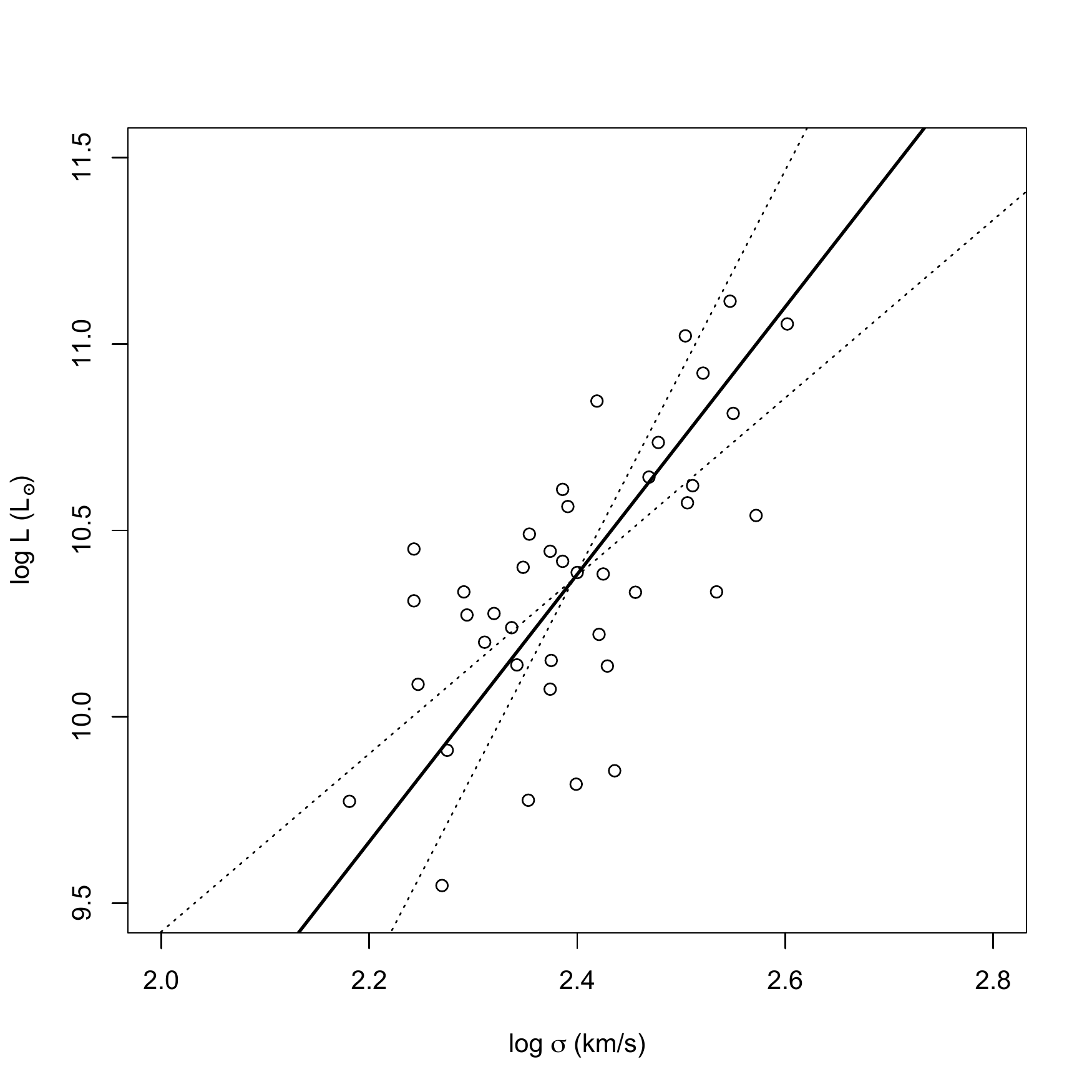}
\caption{Luminosities and velocity dispersions obtained from Schechter's \citeyearpar{schechter1980} measurements of the Faber-Jackson relation with the regression line based on the median of \eqref{posteriorbeta} (bold) and the two ordinary least squares regression lines (dotted).}\label{fjplot}
\end{center}
\end{figure}

The following example illustrates the dilemma posed by estimates that do not possess the same symmetries as the problem statement.  The data are the luminosities $L$ and velocity dispersions $\sigma$ of $n=40$ elliptical galaxies obtained from Schechter's \citeyearpar{schechter1980} measurements of the Faber-Jackson relation, $L \sim \sigma^\beta$, as presented by \citet[][Section 4]{isobe1990}.  As these authors explain, theoretical predictions of $\beta$ range from 2 to 3 to 4.  Figure \ref{fjplot} is a log-log plot of $L$ versus $\sigma$.  The figure strongly suggests a linear relationship, although there is considerable scatter, due primarily to uncharacterized intrinsic processes.

The popular OLS bisector and orthogonal regression estimates of $\beta$ used by astrono\-mers in this and other cosmic distance scale applications are invariant under interchange of the coordinates.  The OLS bisector line bisects the angle between the two OLS regression lines, also shown in Figure \ref{fjplot}, with slopes $b_1=S_{12}/S_{11}=\text{2$\cdot$4}$ $(\text{0$\cdot$4})$ and $b_2=S_{22}/S_{12}=\text{5$\cdot$4}$ $(\text{1$\cdot$0})$.  It has slope $b_{\text{OLSB}}(b_1,b_2)=\tan(\tfrac{1}{2}(\arctan b_1+\arctan b_2))=\text{3$\cdot$4}$ $(\text{0$\cdot$4})$.  The orthogonal regression line minimizes the sum of the squared perpendicular distances to the data.  It has slope $b_{\text{OR}}(b_1,b_2)=B+\sign(S_{12})(B^2+1)^{1/2}=\text{5$\cdot$2}$ $(\text{1$\cdot$0})$, where $B=\tfrac{1}{2}(b_2-1/b_1)$.  Here and in the following the standard errors are estimated from $10^4$ bootstrap replicates.

The OLS bisector and orthogonal regression estimates of $\beta$ are not invariant under scaling of the coordinates.  In theoretical developments of the Faber-Jackson relation, the velocity dispersion generally enters raised to the second power via the kinetic energy, and in general one would expect the analysis of $L \sim (\sigma^k)^{\beta/k}$ to lead to the same estimate of $\beta$ for any $k>0$.  This translates under log transformation to scale invariance.  Defining $b_{\text{OLSB}}^{(k)}/k=b_{\text{OLSB}}(b_1/k,b_2/k)$, the limits $b_{\text{OLSB}}^{(0)}=2/(b_1^{-1}+b_2^{-1})=\text{3$\cdot$3}$ $(\text{0$\cdot$4})$ and $b_{\text{OLSB}}^{(\infty)}=(b_1+b_2)/2=\text{3$\cdot$9}$ $(\text{0$\cdot$5})$.  Similarly, defining $b_{\text{OR}}^{(k)}/k=b_{\text{OR}}(b_1/k,b_2/k)$, the limits $b_{\text{OR}}^{(0)}=b_2=\text{5$\cdot$4}$ $(\text{1$\cdot$0})$ and $b_{\text{OR}}^{(\infty)}=b_1=\text{2$\cdot$4}$ $(\text{0$\cdot$4})$.

\begin{figure}[!ht]
\begin{center}
\includegraphics[width=0.85\textwidth]{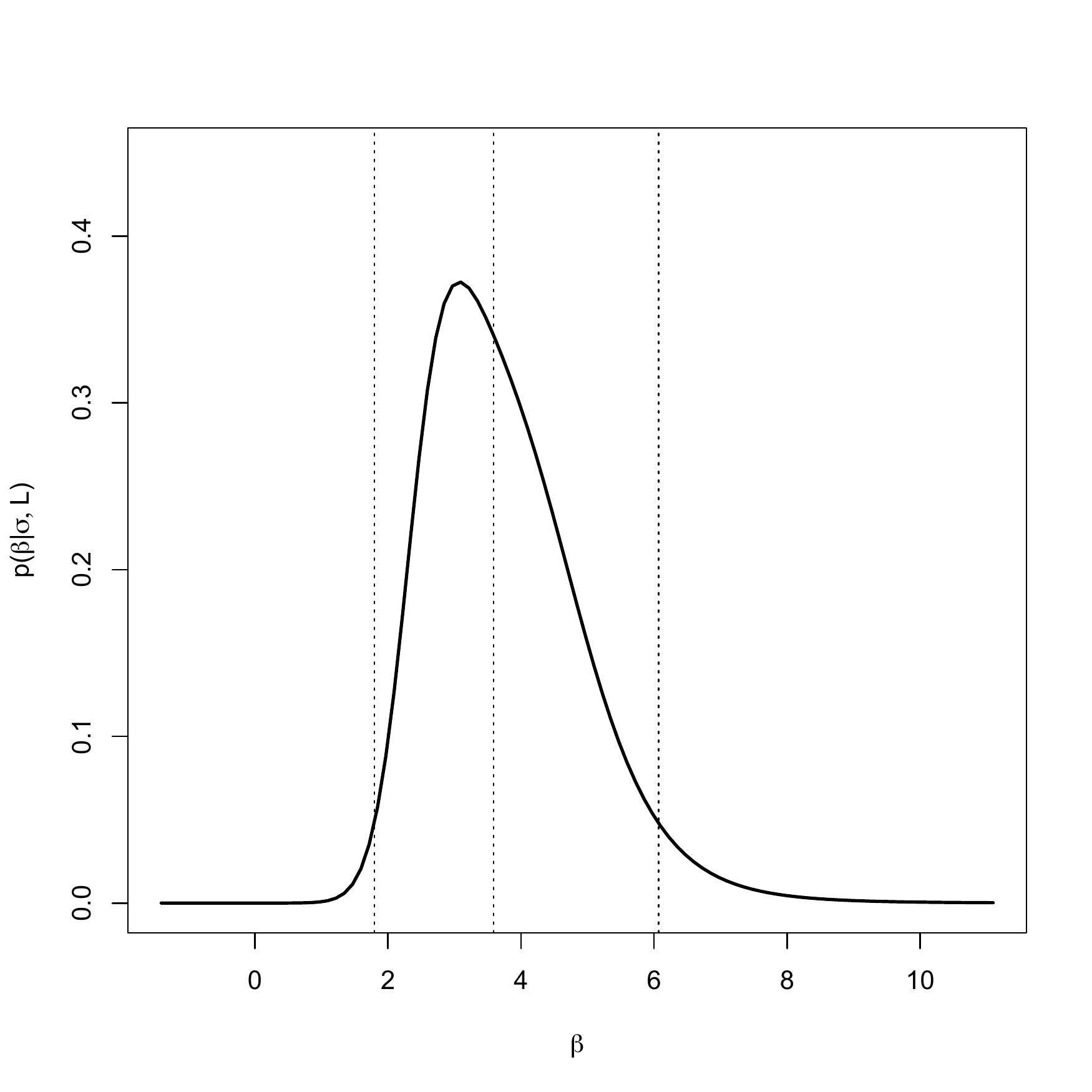}
\caption{Posterior density of $\beta$ in the Faber-Jackson relation $L \sim \sigma^\beta$ with the median and shortest 95\% probability interval.}\label{fjdensity}
\end{center}
\end{figure}

With no criterion of choice, investigators have no firm interval estimate of $\beta$.  \citet{isobe1990} state, ``In cases like these, the astronomer would be wise to calculate [a number of] regressions and be appropriately cautious regarding the confidence of the inferred conclusion.''  The strategy introduced in Section \ref{posterior} inspires a more positive state of affairs.  Figure \ref{fjdensity} illustrates the interchange and scale invariant posterior density \eqref{posteriorbeta}.  The posterior median 3$\cdot$6 favors the theoretical predictions $\beta=3$ and $\beta=4$ more than $\beta=2$, but the shortest 95\% probability interval $(\text{1$\cdot$8},\text{6$\cdot$1})$ confirms that more evidence is needed.

\subsection{Coverage probability}\label{coverage}

Analysis of the Faber-Jackson data of Section~\ref{faberjackson} showed that the location and width of interval estimates around two interchange invariant point estimates of $\beta$ vary with an arbitrary choice of scale.  Another limitation is that, because the sampling density identifies $\Sigma$ but not $\beta$, confidence intervals for point estimates $b(S)$, all functions of sufficient statistics for $\Sigma$, reflect sampling variation about the population values $b(\Sigma)$ but not $\beta$.  The slope cannot be identified with the population value $b(\Sigma)$ unless further restrictions apply.  For instance, the identifying condition for the scale and interchange invariant geometric mean model mentioned in Section~\ref{introduction} is $\sigma_2^2/\sigma_1^2=\Sigma_{22}/\Sigma_{11}$, and the identifying condition for the orthogonal regression model is $\sigma_1^2=\sigma_2^2$.  Coverage of confidence intervals around these point estimates is expected to reach nominal levels when the identifying conditions hold, but otherwise how they will perform is uncertain.

On the other hand, shortest posterior probability intervals calculated from \eqref{posteriorbeta} are marginalized over a distribution of variance parameters.  Considering that the sampling density cannot simultaneously identify the slope and the variance parameters, how the posterior probability intervals will perform when applied to data characterized by a specific configuration of variance parameters is also uncertain.  In this section, I use numerical simulation to study the empirical coverage probability of these intervals.

\begin{table}[htbp]
\begin{center}
\caption{Empirical coverage (\%) of the shortest 90\% posterior probability interval and nominal 90\% confidence intervals of point estimates of $\beta$.}\label{covtbl}\centering\medskip
\begin{tabular}{c @{\extracolsep{0.3cm}} r @{,} l c c c c}
\hline
\hline
\multicolumn{3}{l}{} & Posterior & Geometric & OLS & Orthogonal\\
$n$ & $\sigma_1$&$\sigma_2$ & density & mean & bisector & regression\\
\hline
      & $\text{0$\cdot$05}$&$\text{1$\cdot$00}$ & $\text{86$\cdot$5}$ & $\text{52$\cdot$5}$ & $\text{49$\cdot$8}$ & $\text{77$\cdot$3}$\\
      & $\text{0$\cdot$10}$&$\text{0$\cdot$50}$ & $\text{89$\cdot$9}$ & $\text{79$\cdot$1}$ & $\text{78$\cdot$2}$ & $\text{81$\cdot$1}$\\
20    & $\text{0$\cdot$20}$&$\text{0$\cdot$20}$ & $\text{92$\cdot$8}$ & $\text{84$\cdot$8}$ & $\text{84$\cdot$6}$ & $\text{85$\cdot$3}$\\
      & $\text{0$\cdot$50}$&$\text{0$\cdot$10}$ & $\text{82$\cdot$9}$ & $\text{64$\cdot$0}$ & $\text{63$\cdot$5}$ & $\text{65$\cdot$1}$\\
      & $\text{1$\cdot$00}$&$\text{0$\cdot$05}$ & $\text{72$\cdot$2}$ & $\text{21$\cdot$3}$ & $\text{24$\cdot$4}$ & $\text{20$\cdot$3}$\\
\hline
      & $\text{0$\cdot$05}$&$\text{1$\cdot$00}$ & $\text{80$\cdot$7}$ & $\text{\phantom{0}7$\cdot$8}$ & $\text{\phantom{0}7$\cdot$2}$ & $\text{22$\cdot$4}$\\
      & $\text{0$\cdot$10}$&$\text{0$\cdot$50}$ & $\text{83$\cdot$9}$ & $\text{56$\cdot$3}$ & $\text{55$\cdot$8}$ & $\text{58$\cdot$9}$\\
50    & $\text{0$\cdot$20}$&$\text{0$\cdot$20}$ & $\text{94$\cdot$6}$ & $\text{88$\cdot$0}$ & $\text{88$\cdot$0}$ & $\text{88$\cdot$3}$\\
      & $\text{0$\cdot$50}$&$\text{0$\cdot$10}$ & $\text{75$\cdot$8}$ & $\text{40$\cdot$4}$ & $\text{40$\cdot$4}$ & $\text{41$\cdot$0}$\\
      & $\text{1$\cdot$00}$&$\text{0$\cdot$05}$ & $\text{54$\cdot$4}$ & $\text{\phantom{0}3$\cdot$2}$ & $\text{\phantom{0}3$\cdot$4}$ & $\text{\phantom{0}2$\cdot$8}$\\
\hline
      & $\text{0$\cdot$05}$&$\text{1$\cdot$00}$ & $\text{71$\cdot$6}$ & $\text{\phantom{0}0$\cdot$3}$ & $\text{\phantom{0}0$\cdot$2}$ & $\text{\phantom{0}0$\cdot$9}$\\
      & $\text{0$\cdot$10}$&$\text{0$\cdot$50}$ & $\text{75$\cdot$5}$ & $\text{30$\cdot$2}$ & $\text{30$\cdot$1}$ & $\text{30$\cdot$9}$\\
100   & $\text{0$\cdot$20}$&$\text{0$\cdot$20}$ & $\text{96$\cdot$6}$ & $\text{88$\cdot$0}$ & $\text{88$\cdot$0}$ & $\text{88$\cdot$1}$\\
      & $\text{0$\cdot$50}$&$\text{0$\cdot$10}$ & $\text{71$\cdot$5}$ & $\text{23$\cdot$3}$ & $\text{23$\cdot$2}$ & $\text{22$\cdot$8}$\\
      & $\text{1$\cdot$00}$&$\text{0$\cdot$05}$ & $\text{42$\cdot$1}$ & $\text{\phantom{0}0$\cdot$1}$ & $\text{\phantom{0}0$\cdot$1}$ & $\text{\phantom{0}0$\cdot$0}$\\
\hline
\multicolumn{7}{l}{\footnotesize Coverage based on 1000 random data sets.  Geometric mean, OLS bisector,}\tabularnewline
\multicolumn{7}{l}{\footnotesize and orthogonal regression basic bootstrap confidence intervals estimated}\tabularnewline
\multicolumn{7}{l}{\footnotesize using 999 bootstrap replicates of each data set.  Posterior density from \eqref{posteriorbeta}.}\tabularnewline
\multicolumn{7}{l}{\footnotesize  Geometric mean estimate: $\sign(S_{12})\sqrt{b_1 b_2}$, $b_1=S_{12}/S_{11}$, $b_2=S_{22}/S_{12}$.}\tabularnewline
\multicolumn{7}{l}{\footnotesize OLS bisector estimate: $\tan(\tfrac{1}{2}(\theta_1+\theta_2))$, $\theta_1=\arctan b_1$, $\theta_2=\arctan b_2$.}\tabularnewline
\multicolumn{7}{l}{\footnotesize Orthogonal regression estimate: $B+\sign(S_{12})\sqrt{B^2+1}$, $B=\tfrac{1}{2}(b_2-1/b_1)$.}
\end{tabular}
\end{center}
\end{table}

I considered a number of sample size and measurement error settings, as shown in Table~\ref{covtbl}.  For each setting, I generated 1000 data sets by randomly drawing the true values $\{\xi_{1i}\}\sim \Nrml(0,1)$ and the errors $\{u_{1i}\}\sim \Nrml(0,\sigma_1^2)$ and $\{u_{2i}\}\sim \Nrml(0,\sigma_2^2)$, for $i=1,\dotsc,n$.  I then constructed the observations $\{y_{1i}\}$ and $\{y_{2i}\}$ from the model \eqref{observations} and \eqref{line}, with slope $\beta=1$ and intercept $\alpha=0$.  I calculated the shortest 90\% posterior probability interval for $\beta$ using \eqref{posteriorbeta}, and 90\% basic bootstrap confidence intervals using 999 replicates of each data set for the geometric mean, OLS bisector, and orthogonal regression estimates \citep[see, for example,][Chapter 5]{davisonbook}.  Table~\ref{covtbl} shows the percentage of intervals that contained $\beta=1$.

Table~\ref{covtbl} shows that coverage of confidence intervals around the popular interchange-invariant point estimates of $\beta$ reached the nominal level when $\sigma_1=\sigma_2$, the settings in which, for $\beta=1$ and $\tau^2=1$, the identifying conditions held.  Otherwise, coverage fell well short of the nominal level and worsened in larger samples.

The posterior probability intervals exhibited broader coverage accuracy, overcovering in the vicinity of $\sigma_2^2/\sigma_1^2=\beta^2$ and undercovering in more singular regions of the $(\sigma_1,\sigma_2)$ sampling domain.  Further investigation confirmed that the average coverage of the posterior probability interval over the limiting sampling density $p(\sigma_1^2,\sigma_2^2\mid\beta,\tau^2)\propto \abs{\Sigma}^{-3/2}$ matched the nominal level.

Previous studies have reported much better performance of the geometric mean, OLS bisector, and orthogonal regression estimates (\citealp[][Tables 2 and 3]{babu1992}; \citealp[][Table 8]{warton2006}).  In these studies, however, replicate data was generated from known $\Sigma$, not $\beta$.  Correspondingly, performance was measured relative to the identified value $b(\Sigma)$, not $\beta$.  Evaluated this way, the estimates perform well in general and increasingly well in larger samples, in direct contrast to the present findings.

\subsection{Method comparison}\label{fat}

\begin{figure}[!ht]
\begin{center}
\includegraphics[width=0.85\textwidth]{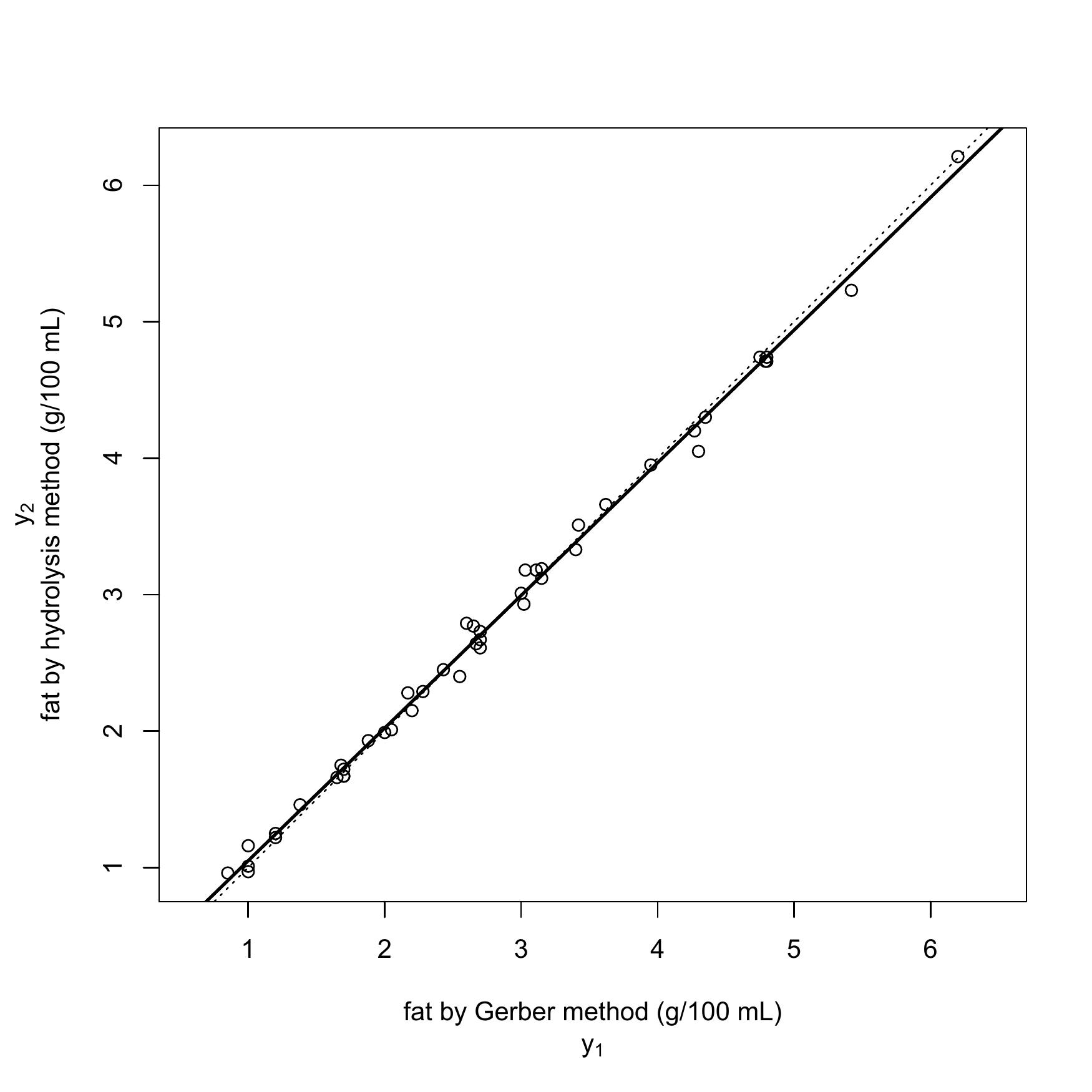}
\caption{Fat content of human milk determined by the standard Gerber method and by enzymic hydrolysis of triglycerides with the regression line based on the median of \eqref{posteriorbeta} (bold) and the line of equality (dotted).}\label{scatter3}
\end{center}
\end{figure}

The final example shows how the posterior density estimate \eqref{posteriorbeta} addresses a limitation of the Bland-Altman approach to method comparison studies \citep{altman1983,bland1986}.  The data are from a study comparing two methods of estimating the fat content of 45 samples of human milk \citep[][Table~3]{bland1999}.  One method ($y_1$) is the standard Gerber method; the other ($y_2$) relies on enzymic hydrolysis of triglycerides.  Figure~\ref{scatter3} shows a scatter plot of the data.

\begin{figure}[!ht]
\begin{center}
\includegraphics[width=0.85\textwidth]{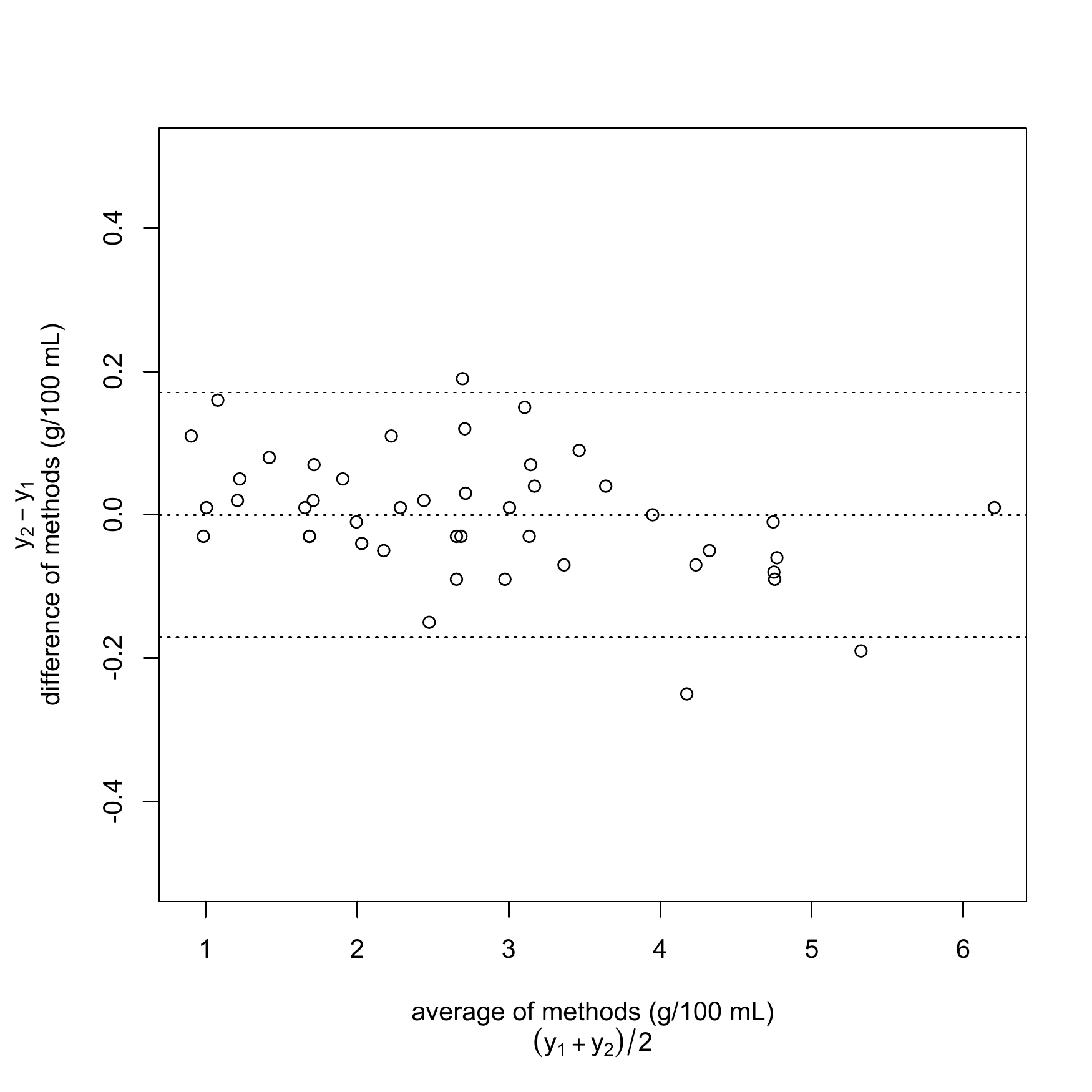}
\caption{Difference versus average of the hydrolysis and Gerber methods of measuring fat content with 95\% limits of agreement.}\label{bias3}
\end{center}
\end{figure}

Standard practice in method comparison relies on the approach described in the highly influential publications of \citet{altman1983,bland1986}.  The centerpiece of the method is the Bland-Altman plot, a plot of the differences against the means of the two methods.  Figure~\ref{bias3} shows a Bland-Altman plot of the fat content data, with horizontal limits of agreement 1.96 standard deviations above and below the mean difference.  The Bland-Altman plot, serving as a type of residuals plot for the identity model, provides a helpful perspective on the data.  In this case, it gives no indication of overall bias, but it suggests a decreasing trend for the differences relative to the means.

The Bland-Altman approach sidesteps the controversial use of regression and correlation \citep{dunn2007} by focusing on the measurements obtained rather than the quantities being measured.  The measurements obtained, however, are affected not only by differences in the quantities being measured but also by the errors.  Indeed, the differences have variance
$\tau^2(\beta-1)^2+(\sigma_1^2+\sigma_2^2)$
under the model \eqref{covtotal}, showing explicitly the contributions of each of these effects.  Importantly, the differences and means have covariance
$1/2\,[\tau^2(\beta^2-1)+(\sigma_2^2-\sigma_1^2)]$,
showing that the measurements may differ systematically even if $\beta=1$, and also that the measurements may agree even if $\beta \neq 1$.

\begin{figure}[!ht]
\begin{center}
\includegraphics[width=0.85\textwidth]{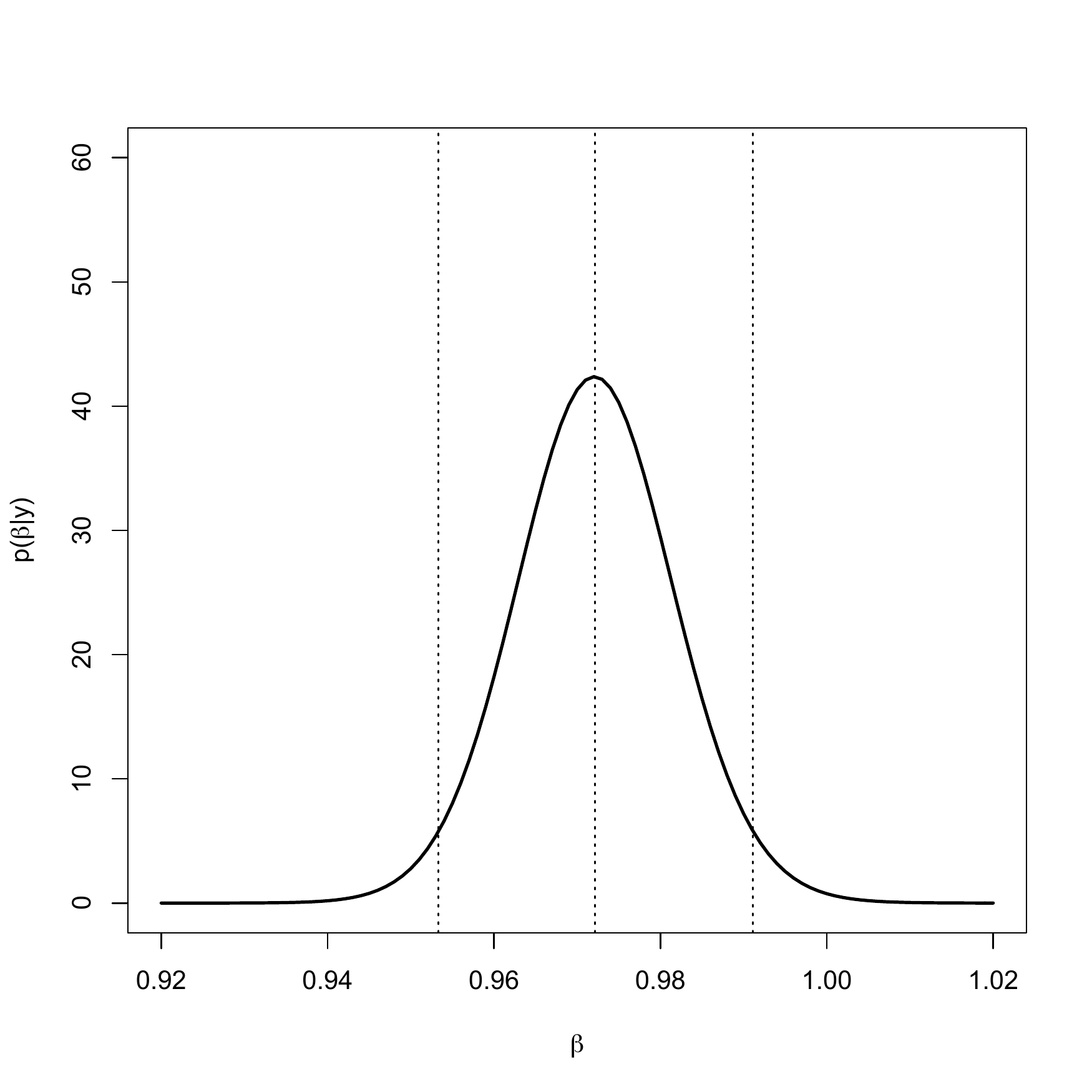}
\caption{Posterior density of the slope $\beta$ for the fat content data with the median and shortest 95\% probability interval.}\label{density3}
\end{center}
\end{figure}

In the present example, the downward trend of Figure~\ref{bias3} hints that $\beta<1$, but  another possibility is that $\sigma_1>\sigma_2$.  The Bland-Altman approach cannot distinguish between these possibilities without specific information on the total errors.  The posterior density shown in Figure~\ref{density3} isolates the relation between the quantities being measured, offering substantial evidence that $\beta<1$, that is, the increments of the quantity being measured by the hydrolysis method are smaller than those of the quantity being measured by the Gerber method.  The shortest 95\% posterior probability interval for $\beta$ is $(\text{0$\cdot$953},\text{0$\cdot$991})$ with median $\text{0$\cdot$972}$.

\section{Discussion}\label{discussion}

Often, the information in a set of observations and knowledge of the sampling density that generated it suffice to identify the parameters of interest.  This is not true in the bivariate normal errors in variables estimation problem with unspecified errors.  The elements of the mean $\mu$ and covariance matrix $\Sigma$ are identified, but the slope $\beta$ is not.  Given $\mu$ and $\Sigma$, the sampling density \eqref{normal} does not vary with $\beta$.  As explained by \citet{poirier1998}, the data are conditionally uninformative for $\beta$.  Fortunately, the nonnegativity conditions \eqref{nonnegativity} modify the domain of $\Sigma$ in a way that depends on $\beta$, so the data are marginally informative for $\beta$.  That is, the data are able to revise prior beliefs, as the examples in Section~\ref{examples} clearly demonstrate.

\citet{reiersol1950} recognized that the nonnegativity conditions \eqref{nonnegativity} restrict the values of $\beta$ given $\Sigma$.  If $\Sigma_{12}>0$, then $\beta$ is restricted to the interval $\Sigma_{12}/\Sigma_{11} \le \beta \le \Sigma_{22}/\Sigma_{12}$.  If $\Sigma_{12}<0$, the inequalities are reversed.  Of course, the true covariance $\Sigma$ is not given, so this fact is not directly useful for inference.  In large samples, however, these bounds should be approximated by the two OLS estimates of $\beta$, an observation \citeauthor{reiersol1950} credited to \citet{frisch1934}, also emphasized by \citet{lindley1968}.

Looking at things the other way around, the nonnegativity conditions \eqref{nonnegativity} restrict the values of $\Sigma$ given $\beta$.  If $\beta>0$, then $\Sigma_{12}$ is restricted to the interval $0 \le \Sigma_{12} \le \min(\beta\Sigma_{11},\allowbreak \Sigma_{22}/\beta)$.  If $\beta<0$, then $\Sigma_{12}$ is restricted to the interval $\max(\beta\Sigma_{11},\Sigma_{22}/\beta) \le \Sigma_{12} \le 0$.  A conditional prior density $p(\Sigma\mid\beta)$ will therefore have a normalization factor that depends on $\beta$, and it will contribute at least this much additional information.  A conditional prior density $p(\Sigma\mid\beta)$ in the inverted Wishart $\W^{-1}(\Psi_0,\nu_0)$ form with $\beta$ dependent normalization factor \eqref{iwishart} will contribute just this much information in the limit of vanishing degrees of freedom; all information that would otherwise be carried by $\Psi_0$ and $\nu_0$ is lost.  Appendix~2 shows that this limit is possible at the very end of the calculation for any proper prior density $p(\beta)$.

The price to pay is that, no matter how large the sample, the prior information carried by the nonnegativity conditions \eqref{nonnegativity} must persist in order to identify $\beta$.  From the perspective of $\beta$ given $\Sigma$, the OLS bounds emphasized by \citet{lindley1968} do not converge.  The sample can never completely overwhelm the joint prior density, and the posterior density cannot concentrate on a single point.

Of course, any information inadvertently incorporated into the prior density will impact the posterior inference on $\beta$ as well, so prior specification cannot be taken lightly.  It is important to consider what is not known as well as what is.  The calculation in Section~\ref{posterior} assumed that prior information on $\beta$ is indifferent to continuous rotations of the coordinate plane, consistent with the rotational symmetry of the sampling density \eqref{normal}.  The calculation went forward by specifying the uniquely rotationally invariant Cauchy prior density \eqref{cauchy} to the scale invariant slope ${\tilde\beta}$.  In contrast, the seemingly benign uniform prior density specified by \citeauthor{zellnerbook} in Section~\ref{zellner} implied a prior density $p(\theta) \propto (\sec\theta)^2$ on the angle $\theta=\arctan{\tilde\beta}$ that introduced a preference for lines approaching the vertical.  Such a preference cannot be desired in all circumstances.

Prior information that is not invariant under continuous rotations of the coordinate plane may still be incorporated in a way that leaves the posterior density invariant under interchange of coordinates.  Interchange invariant prior densities on the scale invariant slope are necessarily of the form $p({\tilde\beta})=1/{\tilde\beta}\,\psi({\tilde\beta},1/{\tilde\beta})$, where $\psi(x_1,x_2)=\psi(x_2,x_1)$ is symmetric with respect to its arguments.  This is a broad class of densities that includes the rotationally invariant Cauchy density and the posterior density \eqref{posteriorbetatilde}.  A convenient way to incorporate prior beliefs on ${\tilde\beta}$ is therefore to use \eqref{posteriorbetatilde} as a prior density, choosing the degree of freedom and correlation parameters $\nu$ and $r$ that reflect these beliefs.  As \eqref{posteriorbetatilde} can be calculated for any $\nu > 1$ and $-1< r < 1$, no prior data are required.

Exact inference from the posterior density \eqref{posteriorbeta} is possible only when the unknown true values and errors are normally distributed.  Whether or not the true values and errors are normally distributed, the normal model has the virtue that the resulting inferences provided by the marginal posterior density depend on these distributions only through the first and second moments of the sampled values, assumed to be finite \citep[see, for example,][Chapter 7]{jaynesbook}.

Section~\ref{problem} described the consequences of a model that assumes either one of the true values of the coordinates is sampled from an improper uniform distribution with infinite variance.  The slope is identified, but the identified value depends on which coordinate is marginalized out, an unacceptable solution to a problem that demands invariance under interchange of coordinates.  By incorporating the assumption that the distributions of true values and errors have finite mean and variance, the normal model leads to a solution with the required symmetries; the drawback is that the sampling density does not identify $\beta$.

Compounding the dilemma, \citet{reiersol1950} showed that the normal model is the only model that does not identify $\beta$.  As \citet{lindley1968} point out, this constitutes an acute sensitivity to distribution.  In any nonnormal model, the limiting posterior density will concentrate; otherwise it will not.  They also point out, however, that more accurate estimation is unlikely without specific information about the true values and errors actually sampled in the data at hand.

Finally, looking forward, it appears to be straightforward to generalize the normal model solution to higher dimensions, starting from the matrix form of the conditional prior density \eqref{iwishart}.  The challenge will be to find a practical expression for the normalization factor \eqref{definition}.  The approach of \citet{klepper1984} may be particularly helpful in this effort.

\appendix
\section*{Appendices}

\subsection*{Appendix 1: The normalization factor}\label{appNormalizing}

The objective is to calculate the normalization factor
\begin{equation}\label{definition}
K(\beta,\nu,\Psi)=\int_R\abs{\Sigma}^{-(\nu+3)/2}\exp{[-\tfrac{1}{2}\tr(\Psi\Sigma^{-1})]\dx{\Sigma}}
\end{equation}
of the density \eqref{iwishart}, where $\Psi$ and $\Sigma$ are positive-definite, symmetric $2 \times 2$ matrices, and $R$ is the region defined by the nonnegativity conditions \eqref{nonnegativity}.  In \eqref{definition} and throughout this appendix, $\Psi$ may depend on $\beta$; for convenience such dependence is suppressed in the notation.  But for the restriction to the integration region $R$, the defining integral in \eqref{definition} could be calculated in much the same way as \citet{fisher1915} first calculated the $2 \times 2$ case of the Wishart density.  As it is, $R$ breaks the defining integral in \eqref{definition} into two separate parts
\begin{equation}\label{partspos}
\begin{split}
K&(\beta,\nu,\Psi)=\\
&\int_{0}^\infty\int_{\beta^2\Sigma_{11}}^\infty\int_{0}^{\beta\Sigma_{11}}\abs{\Sigma}^{-(\nu+3)/2}\exp{[-\tfrac{1}{2}\tr(\Psi\Sigma^{-1})]}\dx{\Sigma_{11}}\dx{\Sigma_{22}}\dx{\Sigma_{12}}\\
&+\int_{0}^\infty\int_{\Sigma_{22}/\beta^2}^\infty\int_{0}^{\Sigma_{22}/\beta}\abs{\Sigma}^{-(\nu+3)/2}\exp{[-\tfrac{1}{2}\tr(\Psi\Sigma^{-1})]}\dx{\Sigma_{22}}\dx{\Sigma_{11}}\dx{\Sigma_{12}},
\end{split}
\end{equation}
for the case $\beta \geq 0$, and
\begin{equation}\label{partsneg}
\begin{split}
K&(\beta,\nu,\Psi)=\\
&\int_{0}^\infty\int_{\beta^2\Sigma_{11}}^\infty\int_{\beta\Sigma_{11}}^{0}\abs{\Sigma}^{-(\nu+3)/2}\exp{[-\tfrac{1}{2}\tr(\Psi\Sigma^{-1})]}\dx{\Sigma_{11}}\dx{\Sigma_{22}}\dx{\Sigma_{12}}\\
&+\int_{0}^\infty\int_{\Sigma_{22}/\beta^2}^\infty\int_{\Sigma_{22}/\beta}^{0}\abs{\Sigma}^{-(\nu+3)/2}\exp{[-\tfrac{1}{2}\tr(\Psi\Sigma^{-1})]}\dx{\Sigma_{22}}\dx{\Sigma_{11}}\dx{\Sigma_{12}},
\end{split}
\end{equation}
for the case $\beta < 0$.

Introducing the reparameterization $(\Psi_{11},\Psi_{22},\Psi_{12}) \rightarrow (\abs{\Psi},r,l)$ of $\Psi$, where
\begin{equation}\label{rl}
\begin{split}
r&=\Psi_{12}/(\Psi_{11}\Psi_{22})^{1/2},\\
l&=(\Psi_{22}/\Psi_{11})^{1/2},
\end{split}
\end{equation}
\eqref{partspos} and \eqref{partsneg} have the form
\begin{equation}\label{combined}
K(\beta,\nu,\Psi)=K_1(\abs{\beta}/l,\nu,\abs{\Psi},r \sign(\beta))+K_1(l/\abs{\beta},\nu,\abs{\Psi},r \sign(\beta)),
\end{equation}
where
\begin{equation}\label{K1sigma}
\begin{split}
K_1({\tilde\beta} &,\nu,\abs{\Psi},r)=\int_{0}^\infty\int_{{\tilde\beta}^2\Sigma_{11}}^\infty\int_{0}^{{\tilde\beta}\Sigma_{11}}\abs{\Sigma}^{-(\nu+3)/2}\\
&\times\exp{[-\tfrac{1}{2}\tfrac{\abs{\Psi}^{1/2}}{\abs{\Sigma}\sqrt{1-r^2}}(\Sigma_{11}-2 r \Sigma_{12}+\Sigma_{22})]}\dx{\Sigma_{11}}\dx{\Sigma_{22}}\dx{\Sigma_{12}}.
\end{split}
\end{equation}

I apply the technique described by \citet[][Section~4.2]{andersonbook} to the integral \eqref{K1sigma}, changing variables $(\Sigma_{22},\Sigma_{12}) \rightarrow (u,v)$ by
\begin{align}
u&=\Sigma_{12}/\Sigma_{11}, \label{u}\\
v&=\Sigma_{22}-\Sigma_{12}^2/\Sigma_{11}. \label{v}
\end{align}
Changing variables $(\Sigma_{11},v) \rightarrow (s,w)$ in the resulting expression by
\begin{align}
s&=\tfrac{1}{2}R_\Psi/\Sigma_{11}, \label{s}\\
w&=\tfrac{1}{2}R_\Psi[1/\Sigma_{11}+Q_r(u)/v], \label{w}
\end{align}
where $R_\Psi=\abs{\Psi}^{1/2}/\sqrt{1-r^2}$, and the nonnegative quadratic function $Q_r(u)=(u-r)^2+(1-r^2)$ leads to
\begin{equation}\label{K1u}
\begin{split}
K_1({\tilde\beta},\nu,\abs{\Psi},r)&=2^\nu\Gamma(\tfrac{\nu+1}{2})\Gamma(\tfrac{\nu-1}{2})R_\Psi^{-\nu}\\
&\times\int_0^{\tilde\beta}Q_r(u)^{-(\nu+1)/2}I_{z(u,{\tilde\beta},r)}(\tfrac{\nu+1}{2},\tfrac{\nu-1}{2})\dx{u}.
\end{split}
\end{equation}
In \eqref{K1u}, $I_z(a,b)$ is the regularized incomplete beta function (beta cumulative distribution function) \citep[][Section 8.17]{dlmf}, and
\begin{equation}\label{Fu}
z(u,{\tilde\beta},r)=\frac{Q_r(u)}{Q_r(u)+{\tilde\beta}^2-u^2}.
\end{equation}

Finally, introducing the integration variable $t$ by
\begin{equation}\label{t}
t/\sqrt{\nu}=(u-r)/\sqrt{1-r^2}
\end{equation}
puts \eqref{K1u} in the form
\begin{equation}\label{K1}
K_1({\tilde\beta},\nu,\abs{\Psi},r)=H(\nu)\abs{\Psi}^{-\nu/2}I({\tilde\beta},\nu,r),
\end{equation}
where
\begin{equation}\label{H}
H(\nu)=\sqrt{\pi}2^\nu\Gamma(\tfrac{\nu}{2})\Gamma(\tfrac{\nu-1}{2}),
\end{equation}
and $I({\tilde\beta},\nu,r)$ is the following integral over the Student $t$ probability density function $p_t(t;\nu)$ with degrees of freedom $\nu$ and $F$ cumulative distribution function $P_F(F;\nu_1,\nu_2)$ with degrees of freedom $\nu_1$ and $\nu_2$ \citep[][Sections 26.6 and 26.7]{handbook}.
\begin{equation}\label{I}
I({\tilde\beta},\nu,r)=\int_{t_-(\nu,r)}^{t_+({\tilde\beta},\nu,r)} p_t(t;\nu)P_F(F(t,{\tilde\beta},\nu,r);\nu+1,\nu-1)\dx{t},
\end{equation}
where the integration limits are
\begin{align}
t_-(\nu,r)&=-\sqrt{\nu}\,r/\sqrt{1-r^2}, \label{tn}\\
t_+({\tilde\beta},\nu,r)&=\sqrt{\nu}\,({\tilde\beta}-r)/\sqrt{1-r^2}, \label{tp}
\end{align}
and
\begin{equation}\label{F}
F(t,{\tilde\beta},\nu,r)=\frac{\nu-1}{\nu+1}\,\frac{\nu+t^2}{[t_+({\tilde\beta},\nu,r)-t_-(\nu,r)]^2-[ t-t_-(\nu,r)]^2}.
\end{equation}
Substituting \eqref{K1} into \eqref{combined}, the final expression for $K(\beta,\nu,\Psi)$ is
\begin{equation}\label{K}
K(\beta,\nu,\Psi)=H(\nu)\abs{\Psi}^{-\nu/2}J(\beta,\nu,r,l),
\end{equation}
where
\begin{equation}\label{J}
J(\beta,\nu,r,l)=I(\abs{\beta}/l,\nu,r\sign(\beta))+I(l/\abs{\beta},\nu,r\sign(\beta)).
\end{equation}

\subsection*{Appendix 2: The noninformative limit}\label{appNoninformative}

Due to the factor $H(\nu)$ in \eqref{K}, the reduced sampling density \eqref{sampling} cannot be evaluated in the limit $\nu_0 \rightarrow 0$.  However, these factors cancel out of the posterior density \eqref{bayes1}, leaving
\begin{equation}\label{bayes2}
\begin{split}
&p(\beta\mid y)=p(\beta)\frac{\abs{\Psi_0(\beta)+\nu S}^{-(\nu_0+\nu)/2}J(\beta,\nu_0+\nu,\Psi_0(\beta)+\nu S)}{\abs{\Psi_0(\beta)}^{-\nu_0/2}J(\beta,\nu_0,\Psi_0(\beta))}\\
&\times\left[\int_{-\infty}^\infty p(\beta)\frac{\abs{\Psi_0(\beta)+\nu S}^{-(\nu_0+\nu)/2}J(\beta,\nu_0+\nu,\Psi_0(\beta)+\nu S)}{\abs{\Psi_0(\beta)}^{-\nu_0/2}J(\beta,\nu_0,\Psi_0(\beta))}\dx{\beta}\right]^{-1}.
\end{split}
\end{equation}
In \eqref{bayes2}, $J(\beta,\nu,\Psi)$ is shorthand for the function $J(\beta,\nu,r,l)$ of \eqref{J} in the parameterization \eqref{rl} of $\Psi$.

The precision matrix $\Psi_0(\beta)$ of \eqref{precision} is parameterized by
\begin{align}
\abs{\Psi_0(\beta)}&=\nu_0^2\kappa_0^4\beta^2(1-\rho_0^2), \label{det0}\\
r_0(\beta)&=\rho_0\sign(\beta), \label{r0}\\
l_0(\beta)&=\abs{\beta}. \label{l0}
\end{align}
The integration limits \eqref{tn} and \eqref{tp} of the first integral $I(\abs{\beta}/l_0,\nu_0,r_0)$ of $J(\beta,\nu_0,r_0,l_0)$ in \eqref{J} are therefore
\begin{align}
t_-(\nu_0,r_0(\beta))&=-\sqrt{\nu_0}\,\rho_0/(1-\rho_0^2)^{1/2}, \label{tn0}\\
t_+(\abs{\beta}/l_0(\beta),\nu_0,r_0(\beta))&=\sqrt{\nu_0}\,(1-\rho_0)/(1-\rho_0^2)^{1/2}, \label{tp0}
\end{align}
independent of $\beta$ whatever the value of $\rho_0$, as is the function $F(t,\beta/l_0,\nu_0,r_0)$, defined in \eqref{F}.  Consequently, $I(\abs{\beta}/l_0,\nu_0,r_0)$ is independent of $\beta$.  The same reasoning can be applied to the second integral $I(l_0/\abs{\beta},\nu_0,r_0)$, and therefore $J(\beta,\nu_0,\Psi_0(\beta))$ cancels out of the posterior density \eqref{bayes2}.  Furthermore, from \eqref{det0}, the factors involving the determinant $\abs{\Psi_0(\beta)}$ in \eqref{bayes2} are well-behaved in the limit
\begin{equation}\label{detlimit}
\lim_{\nu_0 \rightarrow 0}\abs{\Psi_0(\beta)}^{-\nu_0/2}=\lim_{\nu_0 \rightarrow 0}[\nu_0^2\kappa_0^4\beta^2(1-\rho_0^2)]^{-\nu_0/2}=1,
\end{equation}
while $\Psi_0(\beta)=0$ in the same limit.  The noninformative limit of the posterior density \eqref{bayes2} is therefore
\begin{equation}\label{bayes3}
\lim_{\nu_0 \rightarrow 0}p(\beta\mid y)=\frac{p(\beta)J(\beta,\nu,r,l)}{\int_{-\infty}^\infty p(\beta)J(\beta,\nu,r,l)\dx{\beta}},
\end{equation}
where the sample correlation coefficient $r=S_{12}/(S_{11} S_{22})^{1/2}$, and the ratio of standard deviations $l=(S_{22}/S_{11})^{1/2}$ are sufficient statistics.  From the properties of the $F$ cumulative distribution function, the integrals in $J(\beta,\nu,r,l)$ are well-defined for $n>2$.

\subsection*{Appendix 3: Invariance properties}\label{appInvariance}

\citet{samuelson1942} proved that the geometric mean of the two OLS estimates is the only point estimate of the slope consistent with the following three properties: (1) it must depend on the data only through the sample correlation coefficient and ratio of standard deviations; (2) it must be invariant to interchange of the coordinates; (3) it must be invariant to a scale change of either coordinate.

Consider a posterior density $p(\beta\mid y)=f(\beta,r,l)$ exhibiting property~1.  If this density must also exhibit properties~2 and 3, then
\begin{equation}\label{invariance1}
f(\beta,r,l)=f(1/\beta,r,1/l)/\beta^2,
\end{equation}
and
\begin{equation}\label{invariance2}
f(\beta,r,l)=c f(c\beta,r,c l),
\end{equation}
for any $c>0$.  Simultaneous solutions of \eqref{invariance1} and \eqref{invariance2} are of the form
\begin{equation}\label{solution}
f(\beta,r,l)=g({\tilde\beta},1/{\tilde\beta},r)/({\tilde\beta} l),
\end{equation}
where the scale invariant slope ${\tilde\beta}=\beta/l$, and $g$ is any function symmetric in its first two arguments.  The posterior density \eqref{posteriorbeta} is a particular case of \eqref{solution} with
\begin{equation}\label{symmetric}
g({\tilde\beta},1/{\tilde\beta},r)=\frac{{\tilde\beta} p({\tilde\beta})J({\tilde\beta},\nu,r,1)}{h(r)},
\end{equation}
where $h(r)=\int_{-\infty}^\infty p({\tilde\beta})J({\tilde\beta},\nu,r,1)\dx{\tilde\beta}$, and it is easily verified from \eqref{Jbody} and \eqref{cauchy} that ${\tilde\beta} p({\tilde\beta})$ and $J({\tilde\beta},\nu,r,1)$ are each symmetric with respect to ${\tilde\beta}$ and $1/{\tilde\beta}$.

\subsection*{Appendix 4: Special cases}\label{appSpecial}

In special cases, the posterior density \eqref{posteriorbeta} is available in closed form.  For instance, starting from \eqref{K1u} with $n=4$, $\nu=n-1=3$, it is straightforward to show that
\begin{equation}\label{K1uNu3}
K_1({\tilde\beta},3,\abs{\Psi},r)=2^3\abs{\Psi}^{-3/2}(1-r^2)^{3/2}\frac{\tilde\beta}{1+{\tilde\beta}^2}\frac{1}{{\tilde\beta}^2-2r{\tilde\beta}+1}.
\end{equation}
The normalization factor \eqref{combined} becomes
\begin{equation}\label{combinedNu3}
K(\beta,3,\Psi)=2^3\abs{\Psi}^{-3/2}(1-r^2)^{3/2}\frac{\abs{\tilde\beta}}{{\tilde\beta}^2-2r{\tilde\beta}+1},
\end{equation}
where ${\tilde\beta}=\beta/l$ is the scale invariant slope parameter, using the reparameterization \eqref{rl} of $\Psi$.  Applying the limit $\nu_0\rightarrow 0$ to \eqref{bayes1} as described in Appendix~2, the posterior density $p(\beta\mid y)=p({\tilde\beta}\mid y)/l$, where
\begin{equation}\label{posteriorbetatildeNu3}
p({\tilde\beta}\mid y)=K(r) \frac{1}{1+{\tilde\beta}^2}\, \frac{\abs{\tilde\beta}}{{\tilde\beta}^2-2r{\tilde\beta}+1},
\end{equation}
$l=(S_{22}/S_{11})^{1/2}$ is the ratio of standard deviations, $r=S_{12}/(S_{11} S_{22})^{1/2}$ is the sample correlation coefficient, and
\begin{equation}\label{normalizationNu3}
K(r)=F(1,1;\tfrac{3}{2};r^2)^{-1}=\frac{r\sqrt{1-r^2}}{\arcsin{r}},
\end{equation}
which is continuous at $r=0$, with value $K(0)=1$.  In \eqref{normalizationNu3}, $F(a,b;c;z)$ is the Gauss hypergeometric function \citep[][Section 15.2]{dlmf}.

Similarly, in the case $n=6$,
\begin{equation}\label{posteriorbetatildeNu5}
p({\tilde\beta}\mid y)=K(r) \frac{1}{1+{\tilde\beta}^2}\, \frac{\abs{\tilde\beta}({\tilde\beta}^2-r{\tilde\beta}+1)}{({\tilde\beta}^2-2r{\tilde\beta}+1)^2},
\end{equation}
where
\begin{equation}\label{normalizationNu5}
K(r)=F(2,1;\tfrac{3}{2};r^2)^{-1}.
\end{equation}

\bibliography{leonardBA2011}

\vskip 3ex

\section*{Acknowledgements}
The author thanks the Biostatistics Shared Resource at the Harold C. Simmons Comprehensive Cancer Center at the University of Texas Southwestern Medical Center at Dallas for supporting this work.

\end{document}